\input jytex.tex   % available from hep-th
\typesize=10pt \magnification=1200 \baselineskip17truept
%\baselineskip25truept
\footnotenumstyle{arabic} \hsize=6truein\vsize=8.5truein
%\input castess.lab
%\draft
%\leftmargin=1.25in
%\oddleftmargin=.5in
%\evenleftmargin=1.5in
\sectionnumstyle{blank}
\chapternumstyle{blank}
\chapternum=1
\sectionnum=1
\pagenum=0
%\referencestyle{preordered}
% title style follows

\def\begintitle{\pagenumstyle{blank}\parindent=0pt
\begin{narrow}[0.4in]}
\def\endtitle{\end{narrow}\newpage\pagenumstyle{arabic}}

% exercise style follows

\def\beginexercise{\vskip 20truept\parindent=0pt\begin{narrow}[10
truept]}
\def\endexercise{\vskip 10truept\end{narrow}}

% **************    my jyTeX abbreviations   *****************

\def\eql#1{\eqno\eqnlabel{#1}}
\def\ref{\reference}
\def\peq{\puteqn}
\def\pref{\putref}

\def\mgn{\marginnote}
\def\bex{\begin{exercise}}
\def\eex{\end{exercise}}

% *********************** My definitions ************************

\font\open=msbm10 %scaled\magstep1 % For VAX. Borde p195.

 %scaled\magstep1 % For VAX. Borde p195.
%\font\open=msym10 %scaled\magstep1 % For Arbortxt on PC
%\font\opens=msym8 %scaled\magstep1 % For Arbortxt on PC
  % For Arbortxt on PC, and VAX. Borde p199

%\font\smsb=cmss8
\def\StretchRtArr#1{{\count255=0\loop\relbar\joinrel\advance\count255 by1
\ifnum\count255<#1\repeat\rightarrow}}
\def\StretchLtArr#1{\,{\leftarrow\!\!\count255=0\loop\relbar
\joinrel\advance\count255 by1\ifnum\count255<#1\repeat}}

\def\StretchLRtArr#1{\,{\leftarrow\!\!\count255=0\loop\relbar\joinrel\advance
\count255 by1\ifnum\count255<#1\repeat\rightarrow\,\,}}

\def\mbox#1{{\leavevmode\hbox{#1}}}

\def\hspace#1{{\phantom{\mbox#1}}}
\def\oR{\mbox{\open\char82}}

\def\oZ{\mbox{\open\char90}}
\def\oC{\mbox{\open\char67}}

\def\al{\alpha}
 %in jyTeX
 %in jyTeX
 %in jyTeX
 %in jyTeX
 %in jyTeX
 %in jyTeX
 %in jyTeX
 %in jyTeX
 %in jyTeX
% in jyTeX
% in jyTeX
% in jyTeX
\def\be{\beta}
\def\ga{\gamma}
\def\de{\delta}
\def\Ga{\Gamma}

\def\la{\lambda}

\def\om{\omega}

\def\si{\sigma}
\def\Si{\Sigma}
\def\th{\theta}

\def\ze{\zeta}

\def\det{{\rm det\,}}

\def\zf{$\zeta$--function}
\def\zfs{$\zeta$--functions}

     % Newline

\def\frac#1/#2{\leavevmode\kern.1em
\raise.5ex\hbox{\the\scriptfont0 #1}\kern-.1em/\kern-.15em
\lower.25ex\hbox{\the\scriptfont0 #2}}
\def\sfrac#1/#2{\leavevmode\kern.1em
\raise.5ex\hbox{\the\scriptscriptfont0 #1}\kern-.1em/\kern-.15em
\lower.25ex\hbox{\the\scriptscriptfont0 #2}}

\def\gtorder{\mathrel{\raise.3ex\hbox{$>$}\mkern-14mu
             \lower0.6ex\hbox{$\sim$}}}
\def\ltorder{\mathrel{\raise.3ex\hbox{$<$}\mkern-14mu
             \lower0.6ex\hbox{$\sim$}}}

\def\semidirprod{\rlap{\ss C}\raise1pt\hbox{$\mkern.75mu\times$}}
\def\for{\lower6pt\hbox{$\Big|$}}
\def\fish{\kern-.25em{\phantom{abcde}\over \phantom{abcde}}\kern-.25em}

 %triple
%dot
 %double
%dot
 %double dot
%for small #1

\def\boxit#1{\vbox{\hrule\hbox{\vrule\kern3pt
        \vbox{\kern3pt#1\kern3pt}\kern3pt\vrule}\hrule}}
\def\dalemb#1#2{{\vbox{\hrule height .#2pt
        \hbox{\vrule width.#2pt height#1pt \kern#1pt \vrule
                width.#2pt} \hrule height.#2pt}}}

\def\ol{\overline}
        %double stroke
\def\frac#1#2{{{#1}\over{#2}}}
 %lower covariant deriv.
 %upper covariant deriv.
 %lower covariant deriv semicolon.
    %lower ordinary  deriv.
    %lower ordinary  deriv comma.

\def\noin{\noindent}

      %Connection
    %Connection'
\def\comb#1#2{{\left(#1\atop#2\right)}}

\def\cosec{{\rm cosec\,}}
\def\etc{{\it etc. }}

\def\eg{{\it e.g.}}
\def\ie{{\it i.e. }}
\def\cf{{\it cf }}
\def\pa{\partial}

 %gives average <#1>
 %gives thermal average <<#1>>
   %gives bracket <#1|#2>
   %gives comma bracket <#1,#2>
 %gives round bracket (#1,#2)
 %gives round bracket (#1,|#2)
 %gives big bracket <#1|#2>
  %gives
%matrix element <#1|#2|#3>
  %gives reduced matrix element
%<#1||#2||#3>

\def\3j#1#2#3#4#5#6{\left\lgroup\matrix{#1&#2&#3\cr#4&#5&#6\cr}
\right\rgroup}

\def\man{{\cal M}}

\def\m?{\mgn{?}}
% KK's defs

\def\pa{\partial}

\def\beq{\begin{eqnarray}}
\def\eeq{\end{eqnarray}}

%  *******************  Journal refs **********************

\def\aop#1#2#3{{\it Ann. Phys.} {\bf {#1}} ({#2}) #3}

\def\cmp#1#2#3{{\it Comm. Math. Phys.} {\bf {#1}} ({#2}) #3}
\def\cqg#1#2#3{{\it Class. Quant. Grav.} {\bf {#1}} ({#2}) #3}

\def\ijmp#1#2#3{{\it Int. J. Mod. Phys.} {\bf {#1}} ({#2}) #3}

\def\jmp#1#2#3{{\it J. Math. Phys.} {\bf {#1}} ({#2}) #3}
\def\jpa#1#2#3{{\it J. Phys.} {\bf A{#1}} ({#2}) #3}

\def\np#1#2#3{{\it Nucl. Phys.} {\bf B{#1}} ({#2}) #3}
\def\pl#1#2#3{{\it Phys. Lett.} {\bf {#1}} ({#2}) #3}

\def\prp#1#2#3{{\it Phys. Rep.} {\bf {#1}} ({#2}) #3}
\def\pr#1#2#3{{\it Phys. Rev.} {\bf {#1}} ({#2}) #3}
\def\prA#1#2#3{{\it Phys. Rev.} {\bf A{#1}} ({#2}) #3}

\def\prD#1#2#3{{\it Phys. Rev.} {\bf D{#1}} ({#2}) #3}

\def\rmp#1#2#3{{\it Rev. Mod. Phys.} {\bf {#1}} ({#2}) #3}

\def\zfp#1#2#3{{\it Z. f. Phys.} {\bf {#1}} ({#2}) #3}

\def\cras#1#2#3{{\it Comptes Rend. Acad. Sci. (Paris)} {\bf{#1}} (#2) #3}
\def\prs#1#2#3{{\it Proc. Roy. Soc.} {\bf A{#1}} ({#2}) #3}
\def\pcps#1#2#3{{\it Proc. Camb. Phil. Soc.} {\bf{#1}} ({#2}) #3}
\def\mpcps#1#2#3{{\it Math. Proc. Camb. Phil. Soc.} {\bf{#1}} ({#2}) #3}

\def\amsh#1#2#3{{\it Abh. Math. Sem. Ham.} {\bf {#1}} ({#2}) #3}
\def\am#1#2#3{{\it Acta Mathematica} {\bf {#1}} ({#2}) #3}
\def\aim#1#2#3{{\it Adv. in Math.} {\bf {#1}} ({#2}) #3}
\def\ajm#1#2#3{{\it Am. J. Math.} {\bf {#1}} ({#2}) #3}

\def\aom#1#2#3{{\it Ann. of Math.} {\bf {#1}} ({#2}) #3}
\def\cjm#1#2#3{{\it Can. J. Math.} {\bf {#1}} ({#2}) #3}
\def\bams#1#2#3{{\it Bull.Am.Math.Soc.} {\bf {#1}} ({#2}) #3}

\def\cmh#1#2#3{{\it Comm. Math. Helv.} {\bf {#1}} ({#2}) #3}

\def\dmj#1#2#3{{\it Duke Math. J.} {\bf {#1}} ({#2}) #3}
\def\invm#1#2#3{{\it Invent. Math.} {\bf {#1}} ({#2}) #3}

\def\jdg#1#2#3{{\it J. Diff. Geom.} {\bf {#1}} ({#2}) #3}

\def\joa#1#2#3{{\it J. of Algebra} {\bf {#1}} ({#2}) #3}
\def\jram#1#2#3{{\it J. f. reine u. Angew. Math.} {\bf {#1}} ({#2}) #3}
\def\jims#1#2#3{{\it J. Indian. Math. Soc.} {\bf {#1}} ({#2}) #3}
\def\jlms#1#2#3{{\it J. Lond. Math. Soc.} {\bf {#1}} ({#2}) #3}
\def\jmpa#1#2#3{{\it J. Math. Pures. Appl.} {\bf {#1}} ({#2}) #3}
\def\ma#1#2#3{{\it Math. Ann.} {\bf {#1}} ({#2}) #3}

\def\mz#1#2#3{{\it Math. Zeit.} {\bf {#1}} ({#2}) #3}
\def\ojm#1#2#3{{\it Osaka J.Math.} {\bf {#1}} ({#2}) #3}

\def\plb#1#2#3{{\it Phys. Letts.} {\bf {B#1}} ({#2}) #3}
\def\plms#1#2#3{{\it Proc. Lond. Math. Soc.} {\bf {#1}} ({#2}) #3}
\def\pgma#1#2#3{{\it Proc. Glasgow Math. Ass.} {\bf {#1}} ({#2}) #3}
\def\qjm#1#2#3{{\it Quart. J. Math.} {\bf {#1}} ({#2}) #3}

\def\rmjm#1#2#3{{\it Rocky Mountain J. Math.} {\bf {#1}} ({#2}) #3}

\def\tams#1#2#3{{\it Trans.Am.Math.Soc.} {\bf {#1}} ({#2}) #3}

% *******************   Main text *********************
\begin{title}
\vglue 1truein
%\righttext {MUTP/96/23}
%\righttext{hep-th/96}
\vskip15truept
%\leftline{\today}
%\vskip 30truept
\centertext {\Bigfonts \bf $p$-forms on $d$-spherical tesselations}
\vskip10truept \centertext{\Bigfonts \bf}
 \vskip 20truept
\centertext{J.S.Dowker\footnote{dowker@man.ac.uk}} \vskip 7truept
\centertext{\it Theory Group,} \centertext{\it School of Physics and
Astronomy,} \centertext{\it The University of Manchester,} \centertext{\it
Manchester, England} \vskip40truept
\begin{narrow}
The spectral properties of $p$--forms on the fundamental domains of regular
tesselations of the $d$--dimensional sphere are discussed. The degeneracies
for all ranks, $p$, are organised into a double Poincar\'e series which is
explicitly determined. In the particular case of coexact forms of rank
$(d-1)/2$, for odd $d$, it is shown that the heat--kernel expansion
terminates with the constant term, which equals $(-1)^{p+1}/2$ and that the
boundary terms also vanish, all as expected. As an example of the double
domain construction, it is shown that the degeneracies on the sphere are
given by adding the absolute and relative degeneracies on the hemisphere,
again as anticipated. The eta invariant on S$^3/\Ga$ is computed to be
irrational.

The spectral counting function is calculated and the accumulated degeneracy given exactly. A generalised Weyl-Polya conjecture for $p$-forms is suggested and verified.
\end{narrow}
\vskip 5truept
%\righttext {August 1996}
\vskip 60truept
%\righttext{Typeset in \jyTeX}
\vfil
\end{title}
\pagenum=0
\newpage

\section{\bf 1. Introduction.}

In a recent work, [\pref{Dow40}], I have looked at $p$--forms on tesselations
of the three--sphere. In this follow--up, I expand on the higher--dimensional
aspects of the formalism initiated there. The generating functions are
presented for any $p$--form although I concentrate, for actual \zf\
computations, on (coexact) forms of the middle rank, $p=(d-1)/2$, on the odd
$d$--dimensional factored sphere, S$^d/\Ga$. The reason for this is that the
eigenvalues are perfect squares and the expressions for the spectral objects
can be taken a long way in terms of known quantities.

The deck group, $\Ga$, is the complete symmetry group of a regular polytope in
$n$, $=d+1$, dimensions. In another terminology, it is a real reflection
group. These have all been classified.

This paper should be looked upon as a direct continuation of  [\pref{Dow40}]
and, to avoid repetition, I will use, without derivation, any necessary
equations and results of this reference As there, the analysis is presented
as an example of spectral theory in bounded domains that is easily, and
explicitly, managed {\it via} images. It largely consists of bolting together
already existing pieces of knowledge and taking the expressions a little
further than seems to exist in the literature.

The quantum field theory of anti--symmetric fields has a certain importance,
\eg\ [\pref{SandTy,Obukhov,CandT,Reuter}]. I will consider a selection of
spectral objects, such as heat--kernel expansion coefficients, the Casimir
energy and the eta invariant, as examples.

\section{\bf 2. Spectrum and generating functions}
The coexact eigenvalues of the Hodge de--Rham Laplacian, $d\de+\de d$, on the
$d$--sphere are standard,
  $$
  \la^{CE}(p,l)=\big(l+(d+1)/2\big)^2-\big((d-1)/2-p\big)^2\,,
  \quad l=0,1,\ldots\,\,,
  \eql{ceig2}
  $$
which specialises to
  $$
  \mu^{CE}(p,l)=\big(l+p+1\big)^2\,,
  \quad l=0,1,\ldots\,\,,
  \eql{ceigmo}
  $$
for middle rank forms, if $d$ is odd.

The corresponding degeneracies, $g^{CE}_b(p,l)$, are best encoded in the
generating function defined by,
  $$
  g^{CE}_b(p,\si)\equiv\sum_{l=0}^\infty g^{CE}_b(p,l)\,\si^l\,,
  \eql{gf}
  $$
where the label, $b=r$ or $a$, indicates the conditions satisfied by the
$p$--form on the boundary, $\pa\man$, of the fundamental domain, $\man$, for
the action of $\Ga$ on S$^d$. Absolute (`$a$') conditions arise when the form
is symmetric under this action while relative (`$r$') ones originate from
anti--symmetric behaviour. The relation is the duality one,
  $$
  g^{CE}_b\big(p,\si\big)=g^{CE}_{*b}\big(d-1-p,\si\big)\,,
  \eql{du}
  $$
which is a consequence of the $\oR^n$ duality,
  $$
  h^{CCC}_{*b}\big(n-p,\si\big)=h^{CCC}_b\big(p,\si\big)\,,
  \eql{cccd}
  $$
for the closed--coclosed functions, $h^{CCC}$, and the relation between coexact and closed,
  $$
  g^{CE}_b\big(p,\si\big)=h^{CCC}_b\big(p+1,\si\big)\,.
  \eql{identi}
  $$

An important fact is that forms of the middle rank are self--dual in the
sense that $g^{CE}_b\big((d-1)/2,\si\big)=g^{CE}_{*b}\big((d-1)/2,\si\big)$.

Equations are developed in [\pref{Dow40}] that allow the generating function
to be found in closed form in terms of the degrees, ${\bf
d}=(d_1,d_2,\ldots,d_d)$, that define the polytope (reflection) group, $\Ga$.
The case of $p=1$, $d=3$ was treated in detail there. Now I expose the
general result for {\it any} coexact $p$--form in any dimension $d$,
  $$
  g^{CE}_a\big(p,\si\big)=(-1)^{d+1+p}{\sum_{q=0}^{d-1-p}
  \,(-1)^qe_{d-q}\big(\si^{d_1},\ldots,
  \si^{d_d}\big)\over\si^{p+1}\prod_{i=1}^d\big(1-\si^{d_i}\big)}\,,
  \eql{ggen}
  $$
where the $e_q$ are the elementary symmetric functions.

It is useful to write out the relative generating function from the duality
relation (\peq{du}),
  $$\eqalign{
  g^{CE}_r\big(p,\si\big)&=(-1)^p{\sum_{q=0}^p
  \,(-1)^qe_{d-q}\big(\si^{d_1},\ldots,
  \si^{d_d}\big)\over\si^{d-p}\prod_{i=1}^d\big(1-\si^{d_i}\big)}\,,\cr
  &=(-1)^{p+d}{\sum_{q=d-p}^d
  \,(-1)^qe_q\big(\si^{d_1},\ldots,
  \si^{d_d}\big)\over\si^{d-p}\prod_{i=1}^d\big(1-\si^{d_i}\big)}\,.
  }
  \eql{ggen2}
  $$
I derive these expressions later, while developing the formalism.

As in [\pref{Dow40}], the behaviour under the inversion $\si\to1/\si$ is
important. From (\peq{ggen}),
  $$
  g^{CE}_a\big(p,1/\si\big)=(-1)^{p+1}\si^{p+1}{\sum_{q=0}^{d-1-p}
  \,(-1)^qe_q\big(\si^{d_1},\ldots,
  \si^{d_d}\big)\over\prod_{i=1}^d\big(1-\si^{d_i}\big)}\,,
  \eql{ggen3}
  $$
and combined with (\peq{ggen2}), this gives,
  $$
   T^{CE}_a(p,1/\si)-(-1)^dT^{CE}_r(p,\si)=\si^{p-(d-1)/2}\,(-1)^{p+1}\,,
   \eql{tdu}
  $$
after defining,
  $$
  T^{CE}_b(p,\si)=\si^{(d+1)/2}\,g^{CE}_b(p,\si)\,.
  \eql{teepees}
  $$

For self--dual forms, $d$ is odd and (\peq{tdu}) gives the {\it symmetrical}
part of the `cylinder kernel',
  $$
  T^{CE}_b(p,1/\si)+T^{CE}_b(p,\si)=(-1)^{p+1}\,,\quad (b=a,\,r)\,,
  \eql{sdu}
  $$
which can be employed to advantage when evaluating the \zf\ in the next
section.

\section{\bf 3. Zeta functions and heat--kernels}

A useful spectral organising quantity is the coexact \zf,
  $$
  \ze^{CE}_b(p,s)=\sum_{l=0}^\infty {g^{CE}_b(p,l)\over\la^{CE}(p,l)^s}\,,
  $$
because on a manifold, with or without a boundary, there is the general
decomposition,
  $$
  \ze_b(p,s)=\ze_b^{CE}(p,s)+\ze_b^{CE}(p-1,s)\,,
  \eql{zetot}
  $$
of the total \zf\ for the Hodge--de Rham Laplacian.

Only for middle rank forms can $\ze^{CE}$ be related, using (\peq{ceigmo}),
to the generating functions and, for the remainder of this section, I will
make this simplifying restriction so that $d=2p+1$. The \zf\ is then,
  $$\eqalign{
  \ze^{CE}_a(p,s)&={i\Ga(1-2s)\over2\pi}\int_{C_0}d\tau\,(-\tau)^{2s-1}\,
  T_a^{CE}(p,\tau)\cr
  &={i\Ga(1-2s)\over2\pi}\int_{C_0}d\tau\,(-\tau)^{2s-1}\,
  e^{-\tau(d+1)/2}\,g^{CE}_a(p,\tau)
  }
  \eql{zet0}
  $$
where I have introduced $\tau=-\log \si$ and understand, notationally,
$T^{CE}_a(p,\tau)\equiv T^{CE}_a(p,\si)$ and $g^{CE}_a(p,\tau)\equiv
g^{CE}_a(p,\si)$. (I have written the absolute quantity, but this equals the
relative one.)

Therefore, from (\peq{ggen}),
  $$
  \ze^{CE}_a(p,s)=(-1)^p
  {i\Ga(1-2s)\over2\pi}\int_{C_0}\!\!d\tau\,(-\tau)^{2s-1}
  \sum_{q=0}^{d-1-p}\!\! (-1)^q
 {e_{d-q}\big(e^{-d_1\tau},\ldots,
  e^{-d_d\tau}\big)\over\prod_{i=1}^d\big(1-e^{-d_i\tau}\big)}\,.
  \eql{zet1}
  $$

I now recall the integral representation of the Barnes \zf,
  $$\eqalign{ \zeta_d(s,a|{\bf d})=&{i\Gamma(1-s)\over2\pi}\int_{C_0}
  d\tau {e^{-a\tau}
  (-\tau)^{s-1}\over\prod_{i=1}^d\big(1-e^{-d_i\tau}\big)}\,,\cr
  }
  \eql{barn}
  $$
so that (\peq{zet1}) becomes,
  $$\eqalign{
    \ze^{CE}_a(p,s)&=(-1)^p\,\sum_{q=0}^{d-1-p}
    (-1)^q\sum_{{i_1<i_2<\ldots<i_q}\atop=1}^d
    \ze_d\big(2s,\Si d-d_{i_1}-\ldots-d_{i_q}\,|\,{\bf d}\big)\cr
    &=(-1)^p\,\sum_{q=0}^{d-1-p}
    (-1)^q\sum_{{i_1<i_2<\ldots<i_{d-q}}\atop=1}^d
    \ze_d\big(2s,d_{i_1}+\ldots+d_{i_{d-q}}\,|\,{\bf d}\big)\,\cr
     &=(-1)^{p+d}\,\sum_{q=d-p}^d
    (-1)^q
    \sum_{{i_1<i_2<\ldots<i_q}\atop=1}^d
    \ze_d\big(2s,d_{i_1}+\ldots+d_{i_q}\,|\,{\bf d}\big)\,,\cr
    }
    \eql{zeta1}
  $$
by a simple reordering.

As an example, I calculate the values $\ze_a^{CE}(p,-k/2)$, $k\in\oZ$.
Averaging the first and third lines of (\peq{zeta1}),
  $$\eqalign{
   \ze_a^{CE}(p,-k/2)=(-1)^p{k!\over2(d+k)!\prod d_i}\bigg(&
   \sum_{q=0}^{d-1-p}
    (-1)^{q+k}+\sum_{q=d-p}^{d}
    (-1)^q\bigg)\cr
&\times\,\sum_{{i_1<i_2<\ldots<i_q}\atop=1}^d
   B^{(d)}_{d+k}(d_{i_1}+\ldots+d_{i_q}\,|\,{\bf d})\,.\cr
}
  $$

If $k$ is even, the two sums combine,
  $$\eqalign{
   \ze_a^{CE}(p,-k)&={(-1)^p(2k)!\over2(d+2k)!\prod d_i}
   \sum_{q=0}^d
    (-1)^q\!\!\sum_{{i_1<i_2<\ldots<i_q}\atop=1}^d
   B^{(d)}_{d+2k}(d_{i_1}+\ldots+d_{i_q}\,|\,{\bf d})\,.\cr
}
  \eql{rec2}
  $$
Special interest is attached to the value $k=0$,
  $$\eqalign{
   \ze_a^{CE}(p,0)&={(-1)^p\over 2d!\prod d_i}
   \sum_{q=0}^d(-1)^q\sum_{{i_1<i_2<\ldots<i_q}\atop=1}^d
   B^{(d)}_d(d_{i_1}+\ldots+d_{i_q}\,|\,{\bf d})\,.\cr
}
  \eql{rec3}
  $$

To evaluate these, I note that the Barnes \zf\ satisfies the recursion,
[\pref{Barnesa}],
  $$
  \ze_d(s,a+d_i\,|\,{\bf d})=\ze_d(s,a\,|\,{\bf d})
  -\ze_{d-1}(s,a\,|\,\hat{\bf d}_i)\,,
  $$
whose limiting iteration is, [\pref{Dow9}],
  $$\eqalign{
  \ze_d(s,a+d_1+\ldots+d_d\,|\,{\bf d})&
  -\sum_{*=1}^d\ze_d(s,a+d_1+\ldots+*+\ldots+d_d\,|\,{\bf d})\cr
  &+\sum_{*=1}^d\sum_{*=1}^d\ze_d(s,a+d_1+\ldots
  +*+\ldots+*+\ldots+d_d\,|\,{\bf d})\cr
  &\vdots\cr
  &+(-1)^{d-1}\sum_{i=1}^d\ze_d(s,a+d_i\,|\,{\bf d})
  +(-1)^d\ze_d(s,a\,|\,{\bf d})\cr
  &=a^{-s}\,.
  }
  \eql{urec}
  $$

In the first summation, the star denotes that one of the $d$'s is to be
omitted, in turn. In the second summation, every two different pairs of $d$'s
must be successively omitted, and so on. This is the notation of Barnes
[\pref{Barnesa}]. The star summations and omissions correspond to the more
conventional  ordered summations in \eg\ (\peq{zeta1}).

The iteration, (\peq{urec}), on setting $s$ to specific values or by
extracting the poles of the \zf, leads to identities involving the
generalised Bernoulli polynomials. For example, setting $s$ equal to $-2k$,
it is possible to safely equate $a$ to zero and (\peq{urec}) gives the values
of the expressions (\peq{rec2}) and (\peq{rec3}). I find,
  $$\eqalign{
   \ze_a^{CE}(p,-k)&=0\cr
  \ze_a^{CE}(p,0)&={1\over2}\,(-1)^{p+1}\,,
  }
  \eql{rell}
  $$
which show that the coexact middle form heat--kernel coefficients,
$C^{CE}_{k+d/2}$, $k=1,2,\ldots$, vanish and that the constant term,
$C^{CE}_{d/2}$, equals $(-1)^{p+1}/2$. This derivation is related to, but
independent of, that presented in [\pref{Dow40}]. It does not depend on the
fact that the $C_{d/2}(p)$ coefficient for a general $p$--form vanishes on
the {\it doubled} fundamental domain,
  $$
  C^b_{d/2}(p)+C_{d/2}^{*b}(p)=0\,,\quad\forall p\,.
  \eql{ceeid}
  $$
and could be taken as a proof of this fact.

Furthermore, as mentioned at the end of the previous section, the symmetrical
part of the integrand, (\peq{sdu}), produces (\peq{rell}) immediately,
bypassing explicit use of the recursion formulae which I have given, though,
for completeness.

The heat--kernel expansion terminates with the constant term, which
generalises the result in [\pref{CandH}] on the full sphere. It is almost
obvious that factoring the sphere will not alter this fact. The only question
would be the effect of the fixed points.

Adding the two (equal) constants for absolute and relative conditions gives
the constant for the doubled fundamental domain. In particular, the value,
$(-1)^{p+1}$, holds for the full sphere. This agrees with the known value,
[\pref{CandH,ELV,Dow8}].

The remaining heat--kernel coefficients, $C_{k/2}^{CE}$, follow from the
`positive' poles of the coexact \zf\ at $s=(d-k)/2$, for $k=0,1,\ldots,d-1$,
which themselves result from the known poles of the Barnes \zf\ according to
(\peq{zeta1}),
  $$\eqalign{
   C^{CE}_{k/2}=(-1)^p{\Ga\big((d-k)/2\big)
   \over2k!(d-k-1)!\prod d_i}\bigg(&
   \sum_{q=0}^{d-1-p}
    (-1)^q-\sum_{q=d-p}^{d}
    (-1)^{q+k}\bigg)\cr
&\times\,\sum_{{i_1<i_2<\ldots<i_q}\atop=1}^d
   B^{(d)}_k(d_{i_1}+\ldots+d_{i_q}\,|\,{\bf d})\,.\cr
}
  $$
If $k$ is odd, the summations combine to allow the pole part of (\peq{urec})
to come into use showing that the boundary coefficients, \ie $C_{k/2}^{CE}$
for $k$ odd, are zero, again generalising the result in [\pref{Dow40}]. As
before, this conclusion follows more easily from (\peq{sdu}).

There are no such `topological' simplifications or cancellations for the
other values of the \zf, for example for the coexact Casimir energy,
  $$\eqalign{
  E&={1\over2}\,\ze^{CE}_a(p,-1/2)\cr
  &={(-1)^p\over2(d+1)!\prod d_i}
   \sum_{q=0}^p
    (-1)^q\!\!\sum_{{i_1<i_2<\ldots<i_q}\atop=1}^d
   B^{(d)}_{d+1}\big(d_{i_1}+\ldots+d_{i_q}\,|\,{\bf d}\big)\,,
  }
  $$
 and one is reduced to actual computation.
\section{\bf4. Extensions and elaborations}
Although the expression, (\peq{zet0}), for $\ze_b^{CE}(p,s)$ is valid just
for the middle rank forms on S$^d/\Ga$, it has significance for {\it all} $p$
when the factored sphere is realised as the base of a generalised cone in
$\oR^{d+1}$, [\pref{DandKi}], [\pref{BKD}]. I have also referred to this
construction as a bounded M\"obius corner, [\pref{Dow40,Dowl}]. The
separation of variables into radial and angular introduces a term that
effectively cancels the second part of (\peq{ceig2}). Gilkey, [\pref{Gilkey}]
\S4.7.5, refers to the resulting operator as the {\it normalised} spherical
Laplacian.

In this case the $T$ quantities defined in (\peq{teepees}) are {\it bone
fide} cylinder kernels (without propagation significance) and the
corresponding absolute \zf\ is,
   $$\eqalign{
  \ze^{CE}_a(p,s)&={i\Ga(1-2s)\over2\pi}\int_{C_0}d\tau\,(-\tau)^{2s-1}\,
  T^{CE}_a(p,\tau)\cr
  &=(-1)^p{i\Ga(1-2s)\over2\pi}\int_{C_0}\!\!d\tau\,(-\tau)^{2s-1}
  e^{-\big((d-1)/2-p\big)\,\tau}\,\,\times\cr
  &\hspace{****************}\sum_{q=0}^{d-1-p}\!\!
  (-1)^q {e_{d-q}\big(e^{-d_1\tau},\ldots,
  e^{-d_d\tau}\big)\over\prod_{i=1}^d\big(1-e^{-d_i\tau}\big)}\cr
&=(-1)^d\!\!\sum_{q=p+1}^d(-1)^q \!\! \sum_{{i_1<i_2<\ldots<i_q}\atop=1}^d
\!\!\!
    \ze_d\big(2s,{d-1\over2}-p+d_{i_1}+\ldots+d_{i_q}\,|\,{\bf d}\big)\,,
  }
  \eql{zet4}
  $$
with duality giving the relative $\ze_r(p,s)=\ze_a(d-1-p,s)$,
  $$
   \ze^{CE}_r(p,s)=(-1)^d\!\!\sum_{q=d-p}^d(-1)^q \!\!
   \sum_{{i_1<i_2<\ldots<i_q}\atop=1}^d
\!\!\!
    \ze_d\big(2s,p-{d-1\over2}+d_{i_1}+\ldots+d_{i_q}\,|\,{\bf d}\big)\,.
    \eql{zet5}
  $$

These \zfs\ appear as useful intermediate quantities  but have no independent
dynamical significance. In our work on the ball, Dowker and Kirsten
[\pref{DandKi}], they were referred to as `modified' \zfs. The present
results would allow us to extend the ball calculations to factored bases in a
systematic fashion. For example the computations of the scalar functional
determinants reported in [\pref{Dowl}] could be generalised to $p$--forms.

An expression for the modified \zf\ on the full sphere is  given in equn.
(42) in [\pref{DandKi}]. A related formula can be obtained from our present
results by adding the absolute and relative expressions on the hemisphere,
for which all the degrees are one. Hence from (\peq{zet4}) and (\peq{zet5}),
  $$\eqalign{
  \ze^{CE}_{sphere}(p,s)=(-1)^d\bigg(\sum_{q=p+1}^d(-1)^q&\comb
  dq\,\ze_d\big(2s,{d-1\over2}-p+q\,|\,{\bf d}\big)+\cr
  &\sum_{q=d-p}^d(-1)^q\comb dq\,\ze_d\big(2s,p-{d-1\over2}+q\,|\,{\bf
  d}\big)\bigg)\,,
  }
  $$
which, after some manipulation, is equivalent to the form given in
[\pref{DandKi}].\mgn{CHECK SIGN!}

More generally, adding relative and absolute produces the results for a
`doubled' fundamental domain. Rather than give the full expressions, I will
only look at the consequences of the duality relation, (\peq{tdu}) which
yields,
  $$
  T_{a+r}^{CE}(p,\tau)-(-1)^dT_{a+r}^{CE}(p,-\tau)=2(-1)^{p+d}
  \cosh\big(p-(d-1)/2\big)\tau\,,
  $$
allowing \zf\ values to be easily found as powers of $p-(d-1)/2$. I will not
do this in detail and only remark that these values are independent of the
factoring, $\Ga$.
\section{\bf5. The eta invariant}

An important spectral quantity is the eta invariant, $\eta(0)$, which gives a
measure of the asymmetry of the spectrum. Originally introduced by Atiyah,
Patodi and Singer, [\pref{APS}], as a boundary `correction' to an index, it
has achieved an independent life, and its computation has become a standard
challenge. A number of approaches, simple and sophisticated, are available.
The original one, [\pref{APS}], employs the $G$--index theorem. According to
Donnelly, [\pref{Donnelly}], Millson was the first to evaluate $\eta(0)$ on
lens spaces by direct calculation. A direct spectral computation, in the
particular case of spherical space forms, is also mentioned in [\pref{APS}]
and attributed to Ray. Such a calculation was given, later, by Katase,
[\pref{Katase}], on quotients of the 3--sphere. The analysis is somewhat
involved and the result is just the general angle form given previously. Even
so, I outline my own version below.

The actual numbers for the various homogeneous (fixed point free) quotients
of S$^3$ were computed by Gibbons {\it et al}, [\pref{GPR}], who performed
the group average by summing over the angles that define the elements.
Something similar was done by Seade [\pref{Seade}]. In [\pref{Dow11}], I
offered an algebraic alternative to this rather cumbersome geometric
technique. The eta invariant on spherical space forms has been systematically
investigated by Gilkey [\pref{Gilkey,Gilkey3}].

The eta function, $\eta(s)$, measures the asymmetry of the boundary part of
an operator and, as such, is computed on a closed manifold. As an exercise, I
wish to find it for fundamental domains associated with the quotient
S$^d/\Ga$ and these have a boundary. In fact, things are not quite so bad
because I can work on the doubled fundamental domains resulting from the
restriction to the direct rotational polytope group. However, the domain does
have edges and vertices.

The fundamental domain, $\man$, of the spherical tesselation is the base of
the generalised (metric) cone formed, on $\oR^n$, by the set of reflecting
hyperplanes that define the extended group $\Ga$. As such, it is part of the
boundary of this cone, or M\"obius corner (kaleidoscope), the other part
being the union of the flat sides. Restricting to the rotational subgroup of
$\Ga$ turns the cone into a periodic one whose boundary is just its base, the
doubled fundamental domain, $2\man$. There are, however, singularities of
codimension two, corresponding to the edges of the fundamental domain, and of
codimension three from the vertices, \cf [\pref{DandCh}].

The signature eta function on a $d$--manifold, $N$, (typically a boundary) is
neatly expressed in terms of the middle rank coexact eigenforms, $\phi_l$,
by, [\pref{Cheeger}],
 $$
  \eta(s)=\sum_l \int_N {\phi^*_l\wedge d\phi_l\over \mu^{s+1/2}_l}\,,
 $$
where $\mu_l=\mu^{CE}(p,l)$ of (\peq{ceigmo}), and $p=(d-1)/2$ is odd,
($d=4D-1$). For each label, $l$, there are, possibly, two coexact eigenforms
(`positive' and `negative') that can be chosen\footnote{ The restriction
$d=4D-1$ is necessary for this. I have also allowed for complex eigenforms,
although they can be arranged to be real. I have not distinguished
notationally between the positive and negative types. } to be eigenforms of
$*d$,
  $$
  *d\,\phi_l=\pm \,\om_l\, \phi_l\,,\quad \om_l=\sqrt\mu_l\,,
  $$
and the sum is over both types. The $\om$ spectrum is not generally
symmetric.

Despite my preference for the algebraic method, since all the hard work has
been done on the 3--sphere, I initially use angle summation. Only the
signature eta function will be considered and I now derive, again, the
expression obtained long ago in [\pref{APS}].

In physicist's language, the signature eta function in four dimensions is
just the transverse spin--one spectral asymmetry function on the
three--dimensional boundary, [\pref{Dowqg,Dowit}]. For spherical factors, the
necessary spectral information has been given a number of times before in
various connections, \eg\ [\pref{Dow10,DandJ}], and repeated in
[\pref{Dow11,Dow12}].

In terms of the positive and negative, spin--1 `Hamiltonian' \zfs, on
S$^3/\Ga$,
  $$\eqalign{
  \ze^+(s)&=\sum_{\overline L=1}^\infty {d^+(\ol L)\over(\ol
  L+1) ^s}\cr
  \ze^-(s)&=\sum_{\overline L=3}^\infty {d^-(\ol L)\over(\ol
  L-1) ^s}=\sum_{\overline L=1}^\infty {d^-(\ol L+2)\over(\ol
  L+1) ^s}\,,\cr
  }
  \eql{zetaf}
  $$
the eta function is defined to be,
  $$
  \eta(s)= \ze^+(s)- \ze^-(s)\,.
  $$

The degeneracies, $d^{\pm}$, follow from character theory and are given in
the just cited references. The eta function can be written as the group
average,
  $$
  \eta(s)={1\over|\Ga|}\sum_\ga \eta(\ga,s)\,,
  \eql{ga}
  $$
where, by the algebra detailed in [\pref{Dow12}], the partial eta function,
$\eta(\ga,s)$, is,
  $$
  \eta(\ga,s)=2\sum_{n=1}^\infty{1\over n^s}{\sin\al\sin n\be-\sin\be\sin
  n\al\over\cos\al-\cos\be}\,.
  \eql{ga2}
  $$
The sum over $\ga$ in (\peq{ga}) is a sum over the angles, $\al$ and $\be$.

The important value is $\eta(0)$ and substitution into (\peq{ga2}) shows that
there are no problems with the fixed points. As usual, the identity element
gives zero as do the other special values, $\al=0$, $\be\ne0$, which
correspond to the fixing of a 2--flat in the ambient $\oR^4$. The summation
over $n$ is then trivially performed using $2\sum_{n=1}^\infty \sin n\th=\cot
\th/2$ giving,
  $$
  \eta(0)=-{1\over|\Ga|}\sum_{\al\ne0;\,\be} \cot\al/2\,\cot\be/2\,,
  \eql{etav}
  $$
which is, apart from the summation restriction, the standard formula
\footnote{ This derivation is, no doubt, the same as those of Millson and Ray
mentioned earlier.}.

The values of $\al$ and $\be$ corresponding to the elements of the several
polytope groups can now be inserted and the group average performed using the
class decompositions given in [\pref{Dow40}] which were taken from Hurley
[\pref{Hurley}] and Chang [\pref{Chang}]. I find the values for the doubled
fundamental domain,
  $$\eqalign{
  \eta(0)&=-{2\over5\sqrt5}\,\,,\qquad\qquad\qquad \{3^3\} \cr
  &=-{5\over16}\,,\qquad\qquad\qquad\quad \{3^24\}\cr
  &=-{29\over48}\,,\qquad\qquad\qquad\quad\{343\}\cr
  &=-{2341\over5400}-{118\over75\sqrt5}\,, \qquad\,\{3^25\}\,.
  }
  $$
The novelty is the presence of the surd in two cases and the simple fractions
in the others. By contrast, in the evaluation of the Casimir energy all
irrationalities cancel.

The eta function on a lune is, almost trivially, zero because all group
elements fix a 2--flat `axis' of rotation, and these contribute nothing.

Without going through the mode analysis, it is reasonable that the Dirac eta
invariant will be given by the standard, basic expression as given, in Hanson
and R\"omer, [\pref{HandR}], \eg, see [\pref{EGH}],
  $$
  \eta_S(0)=-{1\over2|\Ga|}\sum_{\al\ne0;\,\be} \cosec\al/2\,\,\cosec\be/2\,.
  \eql{etas}
  $$
Numerical evaluation yields the following values,
  $$\eqalign{
  \eta_S(0)&=-{1\over5\sqrt5}\,\,,\,\,\qquad\qquad\qquad\quad \{3^3\} \cr
  &=-{89\over768}-{9\over32\sqrt2}\,,\qquad\qquad \{3^24\}\cr
  &=-{1867\over1728}-{9\over8 \sqrt2}\,,\qquad\qquad\{343\}\cr
  &=-{37291\over7200}+{277\over75}{1\over\sqrt5}\,, \qquad\,\{3^25\}\,.
  }
  $$

The presence of the surds implies that it is not possible to find alternative
expressions for the eta invariant purely in terms of the degrees, $d_i$, as
it is for {\it homogeneous}, fixed point free quotients, [\pref{Dow11}], or
for the Casimir energy.

Incidentally, this conclusion seems not in agreement with the work of
Degeratu, [\pref{Degeratu}], which relates the coefficients in the Laurent
expansion of the Molien (Poincar\'e) series directly to the (Dirac) eta
invariants associated with the boundaries, S$^3/\Ga_i$, of the orbifolds,
$\oC^2/\Ga_i$ where the $\Ga_i$ are subgroups of $\Ga$.

Cheeger, [\pref{Cheeger}], discusses the eta invariant on a generalised cone.
For the standard situation of a smooth manifold, $M$, a generalised cone is
attached to the boundary, $N$, converting $M$ into $X$, a compact space, with
a conical singularity, on which the index can be calculated. In this way,
Cheeger shows that the standard Atiyah--Patodi--Singer formula for Sig($M$)
follows from spectral analysis on the cone, the boundary $\eta(0)$ arising
now from the effect of the cone apex. Also mentioned is the non--standard eta
function on manifolds with boundaries or with conical points and the
possibility that it might be irrational is raised. My computation seems to
confirm this and I leave it at this point.
\section{\bf6. Developing the formalism -- the double Poincar\'e series}

In this section I present a derivation of my basic formulae, (\peq{ggen}) and
(\peq{ggen2}), from the recursions for the various generating functions given
in [\pref{Dow40}] which are defined as sums over the mode label, as in
(\peq{gf}). Although not necessary, I will do this using double generating
functions obtained from the previous ones by summing also over the form rank,
$p$. Such double series are used by Ray on spheres, [\pref{Ray}], but my
approach is different in detail and, in fact, refers to the expressions {\it
after} the group average. They allow for a compressed treatment.

I start with the degeneracy of harmonic polynomial forms on $\oR^n$ and
define the double Poincar\'e series, a finite, `fermionic' polynomial in $z$,
  $$
  h_b(z,\si)=\sum_{p=0}^n h(p,\si)\,z^p={1-\si^2\over|\Ga|}
  \sum_A{\det(1+zA)\over\det(1-\si A)}\,\chi^*(A)\,.
  \eql{dhpf}
  $$
It is possible to think of $z$ as a fugacity.

Using standard identities, explicit forms are,
  $$
  h_a(z,\si)=(1+z\si)\prod_{i=1}^d{1+z\si^{m_i}
  \over1-\si^{d_i}}\,,
  \eql{dhpf2}
  $$
and
  $$\eqalign{
  h_r(z,\si)&=(z+\si)\prod_{i=1}^d{z+\si^{m_i}
  \over1-\si^{d_i}}\cr
  &=z^n\,h_a\big({1\over z},\si\big)\,,
  }
  \eql{dhpf3}
  $$
where the final factor in the numerator has been extracted using $m_n=1$.
Equation (\peq{dhpf}), without the harmonic factor $(1-\si^2)$, is Solomon's
theorem, [\pref{Solomon}], Bourbaki, [\pref{Bourbaki}] p.136, Kane,
[\pref{Kane}] \S22.4,. See also, relatedly, Sturmfels, [\pref{Sturmfels}],
p.37, Exercise 5. The explicit forms are algebraic while the group average is
geometric.

I introduce the notion of polynomial dual, or reciprocal, by the definition,
  $$
  ^*\!f(z)\equiv z^n\,f(1/z)\,,\quad\ie\,\,\, ^*{}^*f(z)=f(z)\,,
  $$
on a polynomial, $f$, of (unwritten) degree $n$, so that, \eg, the relation
(\peq{dhpf3}) reads,
  $$
  h_{*b}(z,\si)=\,^*h_b(z,\si)\,.
  \eql{hd}
  $$
This helps notationally when dealing with the maximum form rank, $p=n$, on
$\oR^n$.

The recursions given in [\pref{Dow40}] transcribe into formulae that can be
solved algebraically. I give the basic steps. For example, the harmonic and
closed harmonic generating functions are related by the recursion,
  $$
  h_b(p,\si)=h_b^C(p,\si)+\si h^C_b(p+1,\si)
  \eql{rec}
  $$
which, on account of $h^C_b(p,\si)=0$ $(p>n)$, becomes, using the same basic
symbols,
  $$
  h_b(z,\si)=\ol h_b^{\,C}(z,\si)+{\si\over z}\, \big(\ol h^{\,C}_b(z,\si)-
  \ol h^{\,C}_b(0,\si)\big)\,,
   \eql{rec1}
  $$
where $\ol h$ is defined by,
  $$
  \ol h_b^{\,C}(z,\si)=h_b^{\,C}(z,\si)-\si^2\de_{ba}\,.
  \eql{hbar}
  $$

In going from (\peq{rec}) to (\peq{rec1}), the term involving $\si^2$ has
been inserted by hand, {\it via} (\peq{hbar}), for the reason mentioned in
[\pref{Dow40}] for adding a term, $\de_{ba}\de_{p0}\,\de_{l2}$, to the
solution of the recursion for $h^C_b(p,l)$. It is needed to ensure the
required zero mode end point value, $h_b^C(0,\si)=\de_{ba}$, which is not
covered by the exact sequence that provides the recursion.

It is useful at this point to list some of the end point values,
  $$\eqalign{
  &h_b^C(0,\si)=\de_{ba},\,\quad
 ^*\!h^{C}_b(0,\si)=^*\!\!h_b(0,\si),\,\quad h^{CCC}_b(0,\si)=\de_{ba}\,,\cr
 &^*\!h^{CCC}_b(0,\si)=\de_{br},\,\quad ^*\!h^{CC}_b(0,\si)=\de_{br},\,
 \quad h^{CC}_b(0,\si)=h_b(0,\si)\,.
  }
  \eql{epvs}
  $$
The identities established in [\pref{Dow40}] also take on an elegant
appearance. For example the supertraces,
  $$
  h_b(-\si,\si)=\de_{ba}(1-\si^2)\,,\quad h^{CC}_b(-\si,\si)=\de_{ba}\,.
  \eql{sids}
  $$
The first equation is actually an easy consequence of (\peq{dhpf}).

Thus, setting $z=-\si$ in (\peq{rec1}), gives,
  $$
  h_b(-\si,\si)=\ol h_b^C(0,\si)=(1-\si^2)\,\de_{ba}\,,
  $$
as an algebraic check.

The solution of (\peq{rec1}) is,
  $$
  h_b^C(z,\si)={z\over z+\si}\,h_b(z,\si)
  +{1+z\si\over1+z/\si}\,\de_{ba}\,.
  \eql{hch}
  $$

For the closed and the closed--coclosed functions there are no `correction
terms' and the recursion formula, which has the same form as (\peq{rec1}),
gives,
  $$
  h_b^{CCC}(z,\si)={z\over z+\si}\,
  \,h^{CC}_b(z,\si)+{\si\over z+\si}\,\de_{ba}\,,
  \eql{clsph}
  $$
using $h^{CCC}_b(0,\si)=\de_{ba}$.

The left--hand side of (\peq{clsph}) is the quantity required as it encodes
the closed degeneracies on the $d$--sphere, [\pref{IandT}]. The inputs are
the explicit forms (\peq{dhpf2}) and (\peq{dhpf3}) which can be used after
relating $h^{CC}$ and $h^C$ by duality on $\oR^n$. This yields the various
equivalent forms,
  $$\eqalign{
  h^{CC}_b(z,\si)&=z^n\,h^C_{*b}\big({1\over z},
  \si\big)=\,^*\!h^{C}_{*b}(z,\si)\cr
  &={z^n\over 1+z\si}\,h_{*b}\big({1\over
  z},\si\big)+{z^n\si( z+\si)\over 1+z\si}\,\de_{*ba}\cr
  &={1\over 1+z\si}\,h_{b}\big(z,\si\big)+{z^n \si(z+\si)
  \over 1+z\si}\,\de_{br}\,.
  }
  \eql{hcc}
  $$
and
  $$\eqalign{
  ^*\!h^{CC}_b(z,\si)&=h^{C}_{*b}(z,\si)
  ={z\over z+\si}\,h_{*b}\big(z,\si\big)+{1+z\si\over 1+z/\si}\,\de_{br}\,,\cr
  &
  }
  \eql{hcc*}
  $$
using (\peq{hch}). These expressions exhibit immediately the end values in
(\peq{epvs}).

Substitution into (\peq{clsph}) produces,
  $$
  h^{CCC}_b(z,\si)={1\over (1+\si/z)(1+\si z)}\,h_b(z,\si)
  +{z^n z\si\over 1+z\si}\,\de_{br}
  +{\si\over z+\si}\,\de_{ba}\,.
  \eql{hcccz}
  $$

For consistency, duality on $\oR^n$ in the form,
  $$
  h^{CCC}_b(z,\si)=\,^{*\!}h^{CCC}_{*b}(z,\si)\,,
  \eql{hdu}
  $$
can easily be checked. Hence it is sufficient to restrict to absolute
conditions and the explicit form (\peq{dhpf2}) results in the final
expression,
  $$
   h^{CCC}_a(z,\si)={z\over z+\si}\prod_{i=1}^d{1+z\si^{m_i}
  \over1-\si^{d_i}}+{\si\over z+\si}\,,
  \eql{hcccz2}
  $$
which is equivalent to (\peq{ggen}), taking into account the relation
(\peq{identi}).

One way of showing this, is to go backwards and derive the recursion
satisfied by the double Poincar\'e series constructed directly from
(\peq{ggen}), virtually paralleling the previous analysis. For simplicity, I
rename $h^{CCC}_a(*,\si)$ as $H_a(*)$.

From (\peq{ggen}) it is straightforward to derive the form rank recursion,
  $$
   \si^{p+1}H_a(p+1)+\si^p\, H_a(p)={e_p(\si^{d_1},\ldots
   ,\si^{d_d})\over\prod_{i=1}^d(1-\si^{d_i})}\,.
   \eql{frec}
  $$

I next note that, because of the end point value $^{*\!}H(0)=0$, the top
limit in the sum defining $H(z)$ can be put at $p=d$ (appropriate for working
on S$^d$) and the construction of the generating function of $z$ allows the
recursion (\peq{frec}) to be rewritten and then solved, much as before, to
give,
  $$
  H_a(\si z)={z\over 1+z}\prod_{i=1}^d{1+z\si^{d_i}
  \over1-\si^{d_i}}+{1\over1+z}\,,
  $$
which is the same as (\peq{hcccz2}) after setting $z\to z/\si$, moreover, one
can see how the exponents, $m_i$, turn into the degrees, $d_i$. In this,
slightly synthetic way, I have justified the forms (\peq{ggen}) and
(\peq{ggen2}). A direct demonstration is possible.

For convenience I give the expression for the absolute coexact double
generating function that follows from (\peq{hcccz2}) and the relation
(\peq{identi}),
  $$
   g^{CE}_a(z,\si)={1\over z+\si}\bigg[\prod_{i=1}^d{1+z\si^{m_i}
  \over1-\si^{d_i}}-1\bigg]\,,
  \eql{gcea}
  $$
a neat encapsulation of my results. For completeness, the relative version
reads \footnote{ When checking the various relations it is necessary to note
that, algebraically, $g^{CE}_a(p,\si)=1$ for $p=-1$.},
  $$\eqalign{
   g^{CE}_r(z,\si)&={1\over 1+z\si}\bigg[\prod_{i=1}^d{z+\si^{m_i}
  \over1-\si^{d_i}}+z^n\si\bigg]\cr
  &=z^{d-1}\bigg(g^{CE}_a\big({1\over z},\si\big)+z\bigg)\,.
  }
  \eql{gcer}
  $$
I remark again that the $g^{CE}$ are polynomials in $z$.

  On the $d$--hemisphere,
  $$
   g^{CE}_a(z,\si)\bigg|_{hemisphere}={1\over
   z+\si}\bigg[\bigg({1+z\over1-\si}\bigg)^d-1\bigg]\,,
   \eql{cehsa}
  $$
and
  $$
   z^d g^{CE}_r\big({1\over z},\si\big)\bigg|_{hemisphere}={z\over
   z+\si}\bigg[\bigg({1+z\over1-\si}\bigg)^d+{\si\over z}\bigg]\,,
   \eql{cehsr}
  $$
while a supertrace result is,
  $$
   \si \,g^{CE}_a(-\si,\si)+d=\sum_{i=1}^d {1\over1-\si^{d_i}}\,.
  $$

\section{\bf7. The 0--form case}

When $z=0$, \ie $p=0$, the above expressions for $g^{CE}(z,\si)$ do not give
the degeneracies for the {\it general} $\,0$--form. These are, of course,
known, but the most convenient way of including them in the present formalism
is to extend the range of $l$ down to $-1$, which corresponds to a {\it zero
mode}, as is apparent from (\peq{ceig2}). This mode exists only for absolute
(\ie Neumann) conditions and is a uniform function, a polynomial of zero
order. Hence the {\it complete} $0$--form generating function is,
  $$
  h_b(\si)=\de_{ba}+\si\, g^{CE}_b(0,\si)\,,
  $$
which, from (\peq{gcea}) and (\peq{gcer}), yield the known forms,
  $$\eqalign{
  h_a(\si)&=h_N(\si)={1\over\prod_{i=1}^d (1-\si^{d_i})}\cr
  h_r(\si)&=h_D(\si)={\si^{d_0}\over\prod_{i=1}^d (1-\si^{d_i})}\,,
  \quad d_0=\Si_i\, m_i\,.\cr
  }
  $$

\section{\bf8. The double fundamental domain and a check}
Formulae (\peq{gcea}) and (\peq{gcer}) give the generating functions on the
fundamental domain, $\man$. That on the doubled domain, $2\man$, is obtained
by adding these two expressions. This can be checked explicitly for the
hemisphere since the full sphere degeneracies are standard, [\pref{IandT}].
To this end, I split these as in [\pref{DandKi}],
  $$
  g^{CE}(p,l)={(l+d)!\over p!\,(d-p-1)!\,l!}
  \bigg({1\over l+1+p}+{1\over l+d-p}\bigg)\,,\quad l=0,1,\ldots\,\,,
  \eql{stdeg}
  $$
the two parts of which are related by $p\to d-1-p$. Each part corresponds to
a hemisphere contribution. Assuming, for a moment, that this is true, the
combinatorial identity used in [\pref{DandKi}], equn.(41), allows the
eigenvalue generating functions on the hemisphere to be found from
(\peq{stdeg}). Then,
  $$\eqalign{
  g^{CE}_a(p,\si)\bigg|_{hemisphere}&=\sum_{m=p+1}^d
  \comb{m-1}p{1\over(1-\si)^m}\cr
  g^{CE}_r(p,\si)\bigg|_{hemisphere}&=\sum_{m=d-p}^d
  \comb{m-1}{d-p-1}{1\over(1-\si)^m}\,.
  }
  $$
Construction of the double Poincar\'e series turns these into (\peq{cehsa})
and (\peq{cehsr}) by simple algebra, (replace the lower limits by zero),
confirming that the two parts do give the two hemisphere contributions. That
is,
  $$
  g^{CE}(z,\si)\bigg|_{sphere}=g^{CE}_a(z,\si)\bigg|_{hemisphere}
  +g^{CE}_r(z,\si)\bigg|_{hemisphere}\,,
  $$
which is a special case of,
  $$
  g^{CE}(z,\si)\bigg|_{2\man}=g^{CE}_a(z,\si)\bigg|_{\man}
  +g^{CE}_r(z,\si)\bigg|_{\man}\,,
  $$
a $p$--form generalisation of the known scalar result.

As well as the sum of relative and absolute expressions, the difference is
also of interest (see the next section). It is preferable to remove the $z^d$
term in (\peq{gcer}) and define the combination,
 $$
  w(z,\si)=g_r(z,\si)-g_a(z,\si)-z^d
  \eql{diffg}
 $$
which is a polynomial of degree $d-1$ and is anti--reciprocal,
 $$
  {}^*w(z,\si)=-w(z,\si)\,.
  \eql{sr}
 $$
\section{\bf9. The counting function}

The counting function, $N(\la)$, could be considered the basic global
spectral object\footnote{Probably it first appears in Sturm's 1829 treatment
of the roots of the secular equation}. In general terms, let the eigenvalues,
$\la_i$, be ordered linearly $\la_1\le\la_2\le\ldots\le \la_i\le\ldots$ with
$i$ a counting label and degeneracies accounted for by equality. Then a
definition of $N(\la)$ is
  $$
  N(\la)=\sum_{\la_i<\la}1+\sum_{\la_i=\la}{1\over2}\,.
  \eql{cfun1}
  $$
Alternatively, if the spectrum is described by the distinct {\it eigenlevels}
$\la(n)$, $n=0,1,2,\ldots$, with explicit degeneracies,
$g(n)=g\big(\la(n)\big)$, \cf Baltes and Hilf [\pref{BaandH}],
  $$
  N(\la)=\sum_{\la(n)<\la}g(n)+{1\over2}\,g(n)\,\de_{\la,\la(n)}\,.
\eql{cfun2}
  $$
In the case under consideration in this paper, $n$ is $l$, the polynomial
order, and has a dynamical significance. The eigenvalues are functions of
$l$, hence the notation.

The value of $N(\la)$ depends on the particular function $\la(n)$ in the
sense that, for a fixed argument, $\la$, it will vary if the form of $\la(n)$
is changed. One could, for example, add a variable constant to the
eigenvalues. On the other hand, the evaluation of $N(\la)$ at an eigenvalue,
$G(n)\equiv N\big(\la(n)\big)$, is the accumulated degeneracy and does not
depend on the form of $\la(n)$,
  $$
  G(n)=\sum_{n'=0}^n g(n')+{1\over2}g(n)\,.
  $$
$G(n)$ satisfies the recursion,
  $$
  G(n)-G(n-1)={1\over2}\big(g(n)+g(n-1)\big)
  $$
which translates into the relation,
  $$
  G(\si)={1\over2}\bigg({1+\si\over1-\si}\bigg)\,g(\si)\,.
  \eql{cfun3}
  $$
between generating functions, $G(\si)=\sum_{n=0}^\infty G(n)\si^n$ \etc

This equation can be applied to the situation in this paper and, for example,
from (\peq{gcea}) I find,
  $$
  G^{CE}_a(\si)={1\over2}\bigg({1+\si\over1-\si}\bigg)
  {1\over z+\si}\bigg[\prod_{i=1}^d{1+z\si^{m_i}
  \over1-\si^{d_i}}-1\bigg]\,.
  \eql{cfun4}
  $$

It is possible to extract the individual accumulated degeneracies\footnote{
The degeneracies themselves can be extracted but will not be exhibited here.}
by writing,
  $$
  G(l)={1\over2\pi i}\oint_C d\si\,{1\over\si^{l+1}}\,G(\si)
  $$
where the contour $C$ circles the origin. The singularities of $N(\si)$ lie
only at $\si=1$, and it converges as $|\si|\to\infty$, so the contour can be
deformed into one, $C'$, around $\si=1$,
  $$
  G(l)=-{1\over2\pi i}\oint_{C'} d\si\,{1\over\si^{l+1}}\,G(\si)\,,
  $$
and the calculation is one of residues. In principle this gives an explicit
expression for $N(l)$. As an example, I treat the hemisphere using
(\peq{cehsa}) and (\peq{cfun3}),
  $$\eqalign{
  G^{CE}_a(l)&=-{1\over4\pi i}\oint_{C'} d\si\,{1\over\si^{l+1}}
  {1+\si\over(1-\si)(z+\si)}
  \bigg[\bigg({1+z\over1-\si}\bigg)^d-1 \bigg]\cr
  &={(-1)^d\over2d!}\,(1+z)^d{d^d\over d\si^d}
  {1+\si\over(z+\si)\si^{l+1}}\bigg|_{\si=1}-{1\over1+z}\cr
  &=\sum_{r=0}^dA^d_r(z)(l+1)(l+2)
    \ldots(l+r)-{1\over1+z}\,,
  }
  \eql{enel}
  $$
where,
  $$\eqalign{
  A^d_r(z)&={(-1)^{d-r}\over2\, d!}\comb dr(1+z)^d{d^{d-r}\over d\si^{d-r}}{1+\si\over z+\si}
  \bigg|_{\si=1}\cr
  &=(-1)^{d-r}{(1+z)^d\over2(d-r)!r!} {d^{d-r}\over
  d\si^{d-r}}\bigg(1+{1-z\over
  z+\si}\bigg)\bigg|_{\si=1}\cr
  &={1\over d!}(1+z)^{d-1}\,,\quad r=d\cr
  &= {1\over 2r!}(1-z)(1+z)^{r-1}\,,\quad r<d\,.}
  \eql{adr}
  $$

Separating the $r=d$ and $r=0$ terms, (\peq{enel}) reads,
  $$\eqalign{
  G^{CE}_a(l)={(l+1)\ldots(l+d)\over d!}\,&(1+z)^{d-1}-{1\over2}\,+\cr
  &{1-z\over2}\sum_{r=1}^{d-1}{(1+z)^{r-1}\over r!}\,(l+1)(l+2)\ldots(l+r)\,,}
  \eql{geel}
  $$
which is exact and is, correctly, a finite polynomial in $z$.

It should be noted that, because of the final term in the particular
definition of the counting function, (\peq{cfun1}), $G(l)$ is half an
integer.

For a given eigenvalue form, such as (\peq{ceig2}), the counting function,
$N(\la)$ can be determined by the condition that, if $\la(l)\le\la<\la(l+1)$
then $N(\la)=G(l)$.

Particular interest lies in the asymptotic behaviour of $N(\la)$ as
$\la\to\infty$, in relation to Weyl's conjecture. It is tolerably clear that
the leading term will follow from (\peq{geel}) as $l\to \infty$ since, from
(\peq{ceig2}) $\la\sim l^2$ for very large $\la$. The highest power of $l$ is
  $$
  G^{CE}_a(l)\sim {l^d\over d!}\,(1+z)^{d-1}
  $$
which corresponds, after slight manipulation, to Weyl's term
  $$
  N(\la)\sim \comb {d-1}p {|\man|\over(4\pi)^{d/2}}
  {1\over\Ga(1+d/2)}\,\la^{d/2}
  $$
for coexact $p$--forms when $\man$ is the hemisphere. With more algebraic
work, these calculations can be extended to the general tessellation
(\peq{cfun4}).

The point here is that this leading term depends only on the degeneracies and
not on the specific eigenvalues, so long as $\la=l^2+o(l)$. It is possible to
go further and derive a form generalisation of the Weyl--Polya conjecture
which says that, in the scalar case, $N_D(\la)\le N_N(\la)$. Proceeding as
for the scalar case in B\'erard and Besson, [\pref{BandB}], by analogy to the
definition, (\peq{diffg}), I define, the modified difference of accumulated
degeneracies,
  $$
  W(z,\si)={1\over2}{1+\si\over1-\si}\,w(z,\si)
  =G^{CE}_r(z,\si)-G^{CE}_a(z,\si)-{1\over2}{1+\si\over1-\si}\,z^d\,.
  $$
This is an anti--reciprocal polynomial of degree $d-1$ in $z$,
  $$
  {}^*W(z,\si)=-W(z,\si)\,,
  $$
a statement of duality.

Hence, writing the polynomial as,
  $$
  W(z,\si)=w_0+w_1z+\ldots+w_{d-1}z^{d-1}\,,
  $$
one has, $w_i=-w_{d-1-i}$ and, if $d$ is odd, the middle coefficient
$w_{(d-1)/2}$ is zero corresponding to the self--duality of the middle rank
form discussed earlier.

It is then easy to show that $W(z,\si)$ vanishes at $z=1$ and also, for odd
$d$, at $z=-1$.

It is also a fact that the lower half coefficients are all negative and the
upper half all positive \footnote{ Unfortunately I have not yet been able to
prove this, but symbolic manipulation verifies it.},
  $$\eqalign{
  w_i&<0\,,\quad i=0,1,\ldots,d/2-1\cr
  w_i&>0\,,\quad i=d/2,1,\ldots,d-1\,,\cr
  }
  $$
which constitutes our generalisation of the Weyl--Polya conjecture.

\section{\bf10. Conclusion}

The main formal results are the Poincar\'e series (\peq{gcea}), (\peq{gcer}),
for the coexact Laplacian degeneracies on $d$--dimensional fundamental
domains. Some of the specific consequences, such as the termination of the
heat--kernel expansion, are not unexpected. Results involving S$^3$ generally
extend, in some way, to higher odd dimensions.

I note that the heat--kernel expansions are all of the conventional form.
Because of the fixed points, logarithmic terms might have been expected but
it seems that these are hard to generate, [\pref{BandHe}], [\pref{Seeley}].
The \zf\ has the standard meromorphic structure and the implication is that
the image method (the group average) automatically yields the Friedrichs
extension.

The consequences of an irrational eta invariant for the signature have yet to
be determined.
\newpage

\noin{\bf References.} \vskip5truept
\begin{putreferences}
  \ref{DandA}{Dowker,J.S. and Apps, J.S. \cqg{15}{1998}{1121}.}
  \ref{Weil}{Weil,A., {\it Elliptic functions according to Eisenstein
  and Kronecker}, Springer, Berlin, 1976.}
  \ref{Ling}{Ling,C-H. {\it SIAM J.Math.Anal.} {\bf5} (1974) 551.}
  \ref{Ling2}{Ling,C-H. {\it J.Math.Anal.Appl.}(1988).}
 \ref{BMO}{Brevik,I., Milton,K.A. and Odintsov, S.D. \aop{302}{2002}{120}.}
 \ref{KandL}{Kutasov,D. and Larsen,F. {\it JHEP} 0101 (2001) 1.}
 \ref{KPS}{Klemm,D., Petkou,A.C. and Siopsis {\it Entropy
 bounds, monoticity properties and scaling in CFT's}. hep-th/0101076.}
 \ref{DandC}{Dowker,J.S. and Critchley,R. \prD{15}{1976}{1484}.}
 \ref{AandD}{Al'taie, M.B. and Dowker, J.S. \prD{18}{1978}{3557}.}
 \ref{Dow1}{Dowker,J.S. \prD{37}{1988}{558}.}
 \ref{Dow30}{Dowker,J.S. \prD{28}{1983}{3013}.}
 \ref{DandK}{Dowker,J.S. and Kennedy,G. \jpa{}{1978}{}.}
 \ref{Dow2}{Dowker,J.S. \cqg{1}{1984}{359}.}
 \ref{DandKi}{Dowker,J.S. and Kirsten, K. {\it Comm. in Anal. and Geom.
 }{\bf7} (1999) 641.}
 \ref{DandKe}{Dowker,J.S. and Kennedy,G.\jpa{11}{1978}{895}.}
 \ref{Gibbons}{Gibbons,G.W. \pl{60A}{1977}{385}.}
 \ref{Cardy}{Cardy,J.L. \np{366}{1991}{403}.}
 \ref{ChandD}{Chang,P. and Dowker,J.S. \np{395}{1993}{407}.}
 \ref{DandC2}{Dowker,J.S. and Critchley,R. \prD{13}{1976}{224}.}
 \ref{Camporesi}{Camporesi,R. \prp{196}{1990}{1}.}
 \ref{BandM}{Brown,L.S. and Maclay,G.J. \pr{184}{1969}{1272}.}
 \ref{CandD}{Candelas,P. and Dowker,J.S. \prD{19}{1979}{2902}.}
 \ref{Unwin1}{Unwin,S.D. Thesis. University of Manchester. 1979.}
 \ref{Unwin2}{Unwin,S.D. \jpa{13}{1980}{313}.}
 \ref{DandB}{Dowker,J.S.and Banach,R. \jpa{11}{1978}{2255}.}
 \ref{Obhukov}{Obhukov,Yu.N. \pl{109B}{1982}{195}.}
 \ref{Kennedy}{Kennedy,G. \prD{23}{1981}{2884}.}
 \ref{CandT}{Copeland,E. and Toms,D.J. \np {255}{1985}{201}.}
 \ref{ELV}{Elizalde,E., Lygren, M. and Vassilevich,
 D.V. \jmp {37}{1996}{3105}.}
 \ref{Malurkar}{Malurkar,S.L. {\it J.Ind.Math.Soc} {\bf16} (1925/26) 130.}
 \ref{Glaisher}{Glaisher,J.W.L. {\it Messenger of Math.} {\bf18}
(1889) 1.} \ref{Anderson}{Anderson,A. \prD{37}{1988}{536}.}
 \ref{CandA}{Cappelli,A. and D'Appollonio,\pl{487B}{2000}{87}.}
 \ref{Wot}{Wotzasek,C. \jpa{23}{1990}{1627}.}
 \ref{RandT}{Ravndal,F. and Tollesen,D. \prD{40}{1989}{4191}.}
 \ref{SandT}{Santos,F.C. and Tort,A.C. \pl{482B}{2000}{323}.}
 \ref{FandO}{Fukushima,K. and Ohta,K. {\it Physica} {\bf A299} (2001) 455.}
 \ref{GandP}{Gibbons,G.W. and Perry,M. \prs{358}{1978}{467}.}
 \ref{Dow4}{Dowker,J.S..}
  \ref{Rad}{Rademacher,H. {\it Topics in analytic number theory,}
Springer-Verlag,  Berlin,1973.}
  \ref{Halphen}{Halphen,G.-H. {\it Trait\'e des Fonctions Elliptiques},
  Vol 1, Gauthier-Villars, Paris, 1886.}
  \ref{CandW}{Cahn,R.S. and Wolf,J.A. {\it Comm.Mat.Helv.} {\bf 51}
  (1976) 1.}
  \ref{Berndt}{Berndt,B.C. \rmjm{7}{1977}{147}.}
  \ref{Hurwitz}{Hurwitz,A. \ma{18}{1881}{528}.}
  \ref{Hurwitz2}{Hurwitz,A. {\it Mathematische Werke} Vol.I. Basel,
  Birkhauser, 1932.}
  \ref{Berndt2}{Berndt,B.C. \jram{303/304}{1978}{332}.}
  \ref{RandA}{Rao,M.B. and Ayyar,M.V. \jims{15}{1923/24}{150}.}
  \ref{Hardy}{Hardy,G.H. \jlms{3}{1928}{238}.}
  \ref{TandM}{Tannery,J. and Molk,J. {\it Fonctions Elliptiques},
   Gauthier-Villars, Paris, 1893--1902.}
  \ref{schwarz}{Schwarz,H.-A. {\it Formeln und
  Lehrs\"atzen zum Gebrauche..},Springer 1893.(The first edition was 1885.)
  The French translation by Henri Pad\'e is {\it Formules et Propositions
  pour L'Emploi...},Gauthier-Villars, Paris, 1894}
  \ref{Hancock}{Hancock,H. {\it Theory of elliptic functions}, Vol I.
   Wiley, New York 1910.}
  \ref{watson}{Watson,G.N. \jlms{3}{1928}{216}.}
  \ref{MandO}{Magnus,W. and Oberhettinger,F. {\it Formeln und S\"atze},
  Springer-Verlag, Berlin 1948.}
  \ref{Klein}{Klein,F. {\it Lectures on the Icosohedron}
  (Methuen, London, 1913).}
  \ref{AandL}{Appell,P. and Lacour,E. {\it Fonctions Elliptiques},
  Gauthier-Villars,
  Paris, 1897.}
  \ref{HandC}{Hurwitz,A. and Courant,C. {\it Allgemeine Funktionentheorie},
  Springer,
  Berlin, 1922.}
  \ref{WandW}{Whittaker,E.T. and Watson,G.N. {\it Modern analysis},
  Cambridge 1927.}
  \ref{SandC}{Selberg,A. and Chowla,S. \jram{227}{1967}{86}. }
  \ref{zucker}{Zucker,I.J. {\it Math.Proc.Camb.Phil.Soc} {\bf 82 }(1977)
  111.}
  \ref{glasser}{Glasser,M.L. {\it Maths.of Comp.} {\bf 25} (1971) 533.}
  \ref{GandW}{Glasser, M.L. and Wood,V.E. {\it Maths of Comp.} {\bf 25}
  (1971)
  535.}
  \ref{greenhill}{Greenhill,A,G. {\it The Applications of Elliptic
  Functions}, MacMillan, London, 1892.}
  \ref{Weierstrass}{Weierstrass,K. {\it J.f.Mathematik (Crelle)}
{\bf 52} (1856) 346.}
  \ref{Weierstrass2}{Weierstrass,K. {\it Mathematische Werke} Vol.I,p.1,
  Mayer u. M\"uller, Berlin, 1894.}
  \ref{Fricke}{Fricke,R. {\it Die Elliptische Funktionen und Ihre Anwendungen},
    Teubner, Leipzig. 1915, 1922.}
  \ref{Konig}{K\"onigsberger,L. {\it Vorlesungen \"uber die Theorie der
 Elliptischen Funktionen},  \break Teubner, Leipzig, 1874.}
  \ref{Milne}{Milne,S.C. {\it The Ramanujan Journal} {\bf 6} (2002) 7-149.}
  \ref{Schlomilch}{Schl\"omilch,O. {\it Ber. Verh. K. Sachs. Gesell. Wiss.
  Leipzig}  {\bf 29} (1877) 101-105; {\it Compendium der h\"oheren
  Analysis}, Bd.II, 3rd Edn, Vieweg, Brunswick, 1878.}
  \ref{BandB}{Briot,C. and Bouquet,C. {\it Th\`eorie des Fonctions
  Elliptiques}, Gauthier-Villars, Paris, 1875.}
  \ref{Dumont}{Dumont,D. \aim {41}{1981}{1}.}
  \ref{Andre}{Andr\'e,D. {\it Ann.\'Ecole Normale Superior} {\bf 6} (1877)
  265;
  {\it J.Math.Pures et Appl.} {\bf 5} (1878) 31.}
  \ref{Raman}{Ramanujan,S. {\it Trans.Camb.Phil.Soc.} {\bf 22} (1916) 159;
 {\it Collected Papers}, Cambridge, 1927}
  \ref{Weber}{Weber,H.M. {\it Lehrbuch der Algebra} Bd.III, Vieweg,
  Brunswick 190  3.}
  \ref{Weber2}{Weber,H.M. {\it Elliptische Funktionen und algebraische
  Zahlen},
  Vieweg, Brunswick 1891.}
  \ref{ZandR}{Zucker,I.J. and Robertson,M.M.
  {\it Math.Proc.Camb.Phil.Soc} {\bf 95 }(1984) 5.}
  \ref{JandZ1}{Joyce,G.S. and Zucker,I.J.
  {\it Math.Proc.Camb.Phil.Soc} {\bf 109 }(1991) 257.}
  \ref{JandZ2}{Zucker,I.J. and Joyce.G.S.
  {\it Math.Proc.Camb.Phil.Soc} {\bf 131 }(2001) 309.}
  \ref{zucker2}{Zucker,I.J. {\it SIAM J.Math.Anal.} {\bf 10} (1979) 192,}
  \ref{BandZ}{Borwein,J.M. and Zucker,I.J. {\it IMA J.Math.Anal.} {\bf 12}
  (1992) 519.}
  \ref{Cox}{Cox,D.A. {\it Primes of the form $x^2+n\,y^2$}, Wiley,
  New York, 1989.}
  \ref{BandCh}{Berndt,B.C. and Chan,H.H. {\it Mathematika} {\bf42} (1995)
  278.}
  \ref{EandT}{Elizalde,R. and Tort.hep-th/}
  \ref{KandS}{Kiyek,K. and Schmidt,H. {\it Arch.Math.} {\bf 18} (1967) 438.}
  \ref{Oshima}{Oshima,K. \prD{46}{1992}{4765}.}
  \ref{greenhill2}{Greenhill,A.G. \plms{19} {1888} {301}.}
  \ref{Russell}{Russell,R. \plms{19} {1888} {91}.}
  \ref{BandB}{Borwein,J.M. and Borwein,P.B. {\it Pi and the AGM}, Wiley,
  New York, 1998.}
  \ref{Resnikoff}{Resnikoff,H.L. \tams{124}{1966}{334}.}
  \ref{vandp}{Van der Pol, B. {\it Indag.Math.} {\bf18} (1951) 261,272.}
  \ref{Rankin}{Rankin,R.A. {\it Modular forms} CUP}
  \ref{Rankin2}{Rankin,R.A. {\it Proc. Roy.Soc. Edin.} {\bf76 A} (1976) 107.}
  \ref{Skoruppa}{Skoruppa,N-P. {\it J.of Number Th.} {\bf43} (1993) 68 .}
  \ref{Down}{Dowker.J.S. \np {104}{2002}{153}.}
  \ref{Eichler}{Eichler,M. \mz {67}{1957}{267}.}
  \ref{Zagier}{Zagier,D. \invm{104}{1991}{449}.}
  \ref{Lang}{Lang,S. {\it Modular Forms}, Springer, Berlin, 1976.}
  \ref{Kosh}{Koshliakov,N.S. {\it Mess.of Math.} {\bf 58} (1928) 1.}
  \ref{BandH}{Bodendiek, R. and Halbritter,U. \amsh{38}{1972}{147}.}
  \ref{Smart}{Smart,L.R., \pgma{14}{1973}{1}.}
  \ref{Grosswald}{Grosswald,E. {\it Acta. Arith.} {\bf 21} (1972) 25.}
  \ref{Kata}{Katayama,K. {\it Acta Arith.} {\bf 22} (1973) 149.}
  \ref{Ogg}{Ogg,A. {\it Modular forms and Dirichlet series} (Benjamin,
  New York,
   1969).}
  \ref{Bol}{Bol,G. \amsh{16}{1949}{1}.}
  \ref{Epstein}{Epstein,P. \ma{56}{1903}{615}.}
  \ref{Petersson}{Petersson.}
  \ref{Serre}{Serre,J-P. {\it A Course in Arithmetic}, Springer,
  New York, 1973.}
  \ref{Schoenberg}{Schoenberg,B., {\it Elliptic Modular Functions},
  Springer, Berlin, 1974.}
  \ref{Apostol}{Apostol,T.M. \dmj {17}{1950}{147}.}
  \ref{Ogg2}{Ogg,A. {\it Lecture Notes in Math.} {\bf 320} (1973) 1.}
  \ref{Knopp}{Knopp,M.I. \dmj {45}{1978}{47}.}
  \ref{Knopp2}{Knopp,M.I. \invm {}{1994}{361}.}
  \ref{LandZ}{Lewis,J. and Zagier,D. \aom{153}{2001}{191}.}
  \ref{DandK1}{Dowker,J.S. and Kirsten,K. {\it Elliptic functions and
  temperature inversion symmetry on spheres} hep-th/.}
  \ref{HandK}{Husseini and Knopp.}
  \ref{Kober}{Kober,H. \mz{39}{1934-5}{609}.}
  \ref{HandL}{Hardy,G.H. and Littlewood, \am{41}{1917}{119}.}
  \ref{Watson}{Watson,G.N. \qjm{2}{1931}{300}.}
  \ref{SandC2}{Chowla,S. and Selberg,A. {\it Proc.Nat.Acad.} {\bf 35}
  (1949) 371.}
  \ref{Landau}{Landau, E. {\it Lehre von der Verteilung der Primzahlen},
  (Teubner, Leipzig, 1909).}
  \ref{Berndt4}{Berndt,B.C. \tams {146}{1969}{323}.}
  \ref{Berndt3}{Berndt,B.C. \tams {}{}{}.}
  \ref{Bochner}{Bochner,S. \aom{53}{1951}{332}.}
  \ref{Weil2}{Weil,A.\ma{168}{1967}{}.}
  \ref{CandN}{Chandrasekharan,K. and Narasimhan,R. \aom{74}{1961}{1}.}
  \ref{Rankin3}{Rankin,R.A. {} {} ().}
  \ref{Berndt6}{Berndt,B.C. {\it Trans.Edin.Math.Soc}.}
  \ref{Elizalde}{Elizalde,E. {\it Ten Physical Applications of Spectral
  Zeta Function Theory}, \break (Springer, Berlin, 1995).}
  \ref{Allen}{Allen,B., Folacci,A. and Gibbons,G.W. \pl{189}{1987}{304}.}
  \ref{Krazer}{Krazer}
  \ref{Elizalde3}{Elizalde,E. {\it J.Comp.and Appl. Math.} {\bf 118}
  (2000) 125.}
  \ref{Elizalde2}{Elizalde,E., Odintsov.S.D, Romeo, A. and Bytsenko,
  A.A and
  Zerbini,S.
  {\it Zeta function regularisation}, (World Scientific, Singapore,
  1994).}
  \ref{Eisenstein}{Eisenstein}
  \ref{Hecke}{Hecke,E. \ma{112}{1936}{664}.}
  \ref{Terras}{Terras,A. {\it Harmonic analysis on Symmetric Spaces} (Springer,
  New York, 1985).}
  \ref{BandG}{Bateman,P.T. and Grosswald,E. {\it Acta Arith.} {\bf 9}
  (1964) 365.}
  \ref{Deuring}{Deuring,M. \aom{38}{1937}{585}.}
  \ref{Guinand}{Guinand.}
  \ref{Guinand2}{Guinand.}
  \ref{Minak}{Minakshisundaram.}
  \ref{Mordell}{Mordell,J. \prs{}{}{}.}
  \ref{GandZ}{Glasser,M.L. and Zucker, {}.}
  \ref{Landau2}{Landau,E. \jram{}{1903}{64}.}
  \ref{Kirsten1}{Kirsten,K. \jmp{35}{1994}{459}.}
  \ref{Sommer}{Sommer,J. {\it Vorlesungen \"uber Zahlentheorie}
  (1907,Teubner,Leipzig).
  French edition 1913 .}
  \ref{Reid}{Reid,L.W. {\it Theory of Algebraic Numbers},
  (1910,MacMillan,New York).}
  \ref{Milnor}{Milnor, J. {\it Is the Universe simply--connected?},
  IAS, Princeton, 1978.}
  \ref{Milnor2}{Milnor, J. \ajm{79}{1957}{623}.}
  \ref{Opechowski}{Opechowski,W. {\it Physica} {\bf 7} (1940) 552.}
  \ref{Bethe}{Bethe, H.A. \zfp{3}{1929}{133}.}
  \ref{LandL}{Landau, L.D. and Lishitz, E.M. {\it Quantum
  Mechanics} (Pergamon Press, London, 1958).}
  \ref{GPR}{Gibbons, G.W., Pope, C. and R\"omer, H., \np{157}{1979}{377}.}
  \ref{Jadhav}{Jadhav,S.P. PhD Thesis, University of Manchester 1990.}
  \ref{DandJ}{Dowker,J.S. and Jadhav, S. \prD{39}{1989}{1196}.}
  \ref{CandM}{Coxeter, H.S.M. and Moser, W.O.J. {\it Generators and
  relations of finite groups} Springer. Berlin. 1957.}
  \ref{Coxeter2}{Coxeter, H.S.M. {\it Regular Complex Polytopes},
   (Cambridge University Press,
  Cambridge, 1975).}
  \ref{Coxeter}{Coxeter, H.S.M. {\it Regular Polytopes}.}
  \ref{Stiefel}{Stiefel, E., J.Research NBS {\bf 48} (1952) 424.}
  \ref{BandS}{Brink and Satchler {\it Angular momentum theory}.}
  %\ref{Racah1}
  \ref{Rose}{Rose}
  \ref{Schwinger}{Schwinger,J.}
  \ref{Bromwich}{Bromwich, T.J.I'A. {\it Infinite Series},
  (Macmillan, 1947).}
  \ref{Ray}{Ray,D.B. \aim{4}{1970}{109}.}
  \ref{Ikeda}{Ikeda,A. {\it Kodai Math.J.} {\bf 18} (1995) 57.}
  \ref{Kennedy}{Kennedy,G. \prD{23}{1981}{2884}.}
  \ref{Ellis}{Ellis,G.F.R. {\it General Relativity} {\bf2} (1971) 7.}
  \ref{Dow8}{Dowker,J.S. \cqg{20}{2003}{L105}.}
  \ref{IandY}{Ikeda, A and Yamamoto, Y. \ojm {16}{1979}{447}.}
  \ref{BandI}{Bander,M. and Itzykson,C. \rmp{18}{1966}{2}.}
  \ref{Schulman}{Schulman, L.S. \pr{176}{1968}{1558}.}
  \ref{Bar1}{B\"ar,C. {\it Arch.d.Math.}{\bf 59} (1992) 65.}
  \ref{Bar2}{B\"ar,C. {\it Geom. and Func. Anal.} {\bf 6} (1996) 899.}
  \ref{Vilenkin}{Vilenkin, N.J. {\it Special functions},
  (Am.Math.Soc., Providence, 1968).}
  \ref{Talman}{Talman, J.D. {\it Special functions} (Benjamin,N.Y.,1968).}
  \ref{Miller}{Miller,W. {\it Symmetry groups and their applications}
  (Wiley, N.Y., 1972).}
  \ref{Dow3}{Dowker,J.S. \cmp{162}{1994}{633}.}
  \ref{Cheeger}{Cheeger, J. \jdg {18}{1983}{575}.}
  \ref{Dow6}{Dowker,J.S. \jmp{30}{1989}{770}.}
  \ref{Dow20}{Dowker,J.S. \jmp{35}{1994}{6076}.}
  \ref{Dow21}{Dowker,J.S. {\it Heat kernels and polytopes} in {\it
   Heat Kernel Techniques and Quantum Gravity}, ed. by S.A.Fulling,
   Discourses in Mathematics and its Applications, No.4, Dept.
   Maths., Texas A\&M University, College Station, Texas, 1995.}
  \ref{Dow9}{Dowker,J.S. \jmp{42}{2001}{1501}.}
  \ref{Dow7}{Dowker,J.S. \jpa{25}{1992}{2641}.}
  \ref{Warner}{Warner.N.P. \prs{383}{1982}{379}.}
  \ref{Wolf}{Wolf, J.A. {\it Spaces of constant curvature},
  (McGraw--Hill,N.Y., 1967).}
  \ref{Meyer}{Meyer,B. \cjm{6}{1954}{135}.}
  \ref{BandB}{B\'erard,P. and Besson,G. {\it Ann. Inst. Four.} {\bf 30}
  (1980) 237.}
  \ref{PandM}{Polya,G. and Meyer,B. \cras{228}{1948}{28}.}
  \ref{Springer}{Springer, T.A. Lecture Notes in Math. vol 585 (Springer,
  Berlin,1977).}
  \ref{SeandT}{Threlfall, H. and Seifert, W. \ma{104}{1930}{1}.}
  \ref{Hopf}{Hopf,H. \ma{95}{1925}{313}. }
  \ref{Dow}{Dowker,J.S. \jpa{5}{1972}{936}.}
  \ref{LLL}{Lehoucq,R., Lachi\'eze-Rey,M. and Luminet, J.--P. {\it
  Astron.Astrophys.} {\bf 313} (1996) 339.}
  \ref{LaandL}{Lachi\'eze-Rey,M. and Luminet, J.--P.
  \prp{254}{1995}{135}.}
  \ref{Schwarzschild}{Schwarzschild, K., {\it Vierteljahrschrift der
  Ast.Ges.} {\bf 35} (1900) 337.}
  \ref{Starkman}{Starkman,G.D. \cqg{15}{1998}{2529}.}
  \ref{LWUGL}{Lehoucq,R., Weeks,J.R., Uzan,J.P., Gausman, E. and
  Luminet, J.--P. \cqg{19}{2002}{4683}.}
  \ref{Dow10}{Dowker,J.S. \prD{28}{1983}{3013}.}
  \ref{BandD}{Banach, R. and Dowker, J.S. \jpa{12}{1979}{2527}.}
  \ref{Jadhav2}{Jadhav,S. \prD{43}{1991}{2656}.}
  \ref{Gilkey}{Gilkey,P.B. {\it Invariance theory,the heat equation and
  the Atiyah--Singer Index theorem} (CRC Press, Boca Raton, 1994).}
  \ref{BandY}{Berndt,B.C. and Yeap,B.P. {\it Adv. Appl. Math.}
  {\bf29} (2002) 358.}
  \ref{HandR}{Hanson,A.J. and R\"omer,H. \pl{80B}{1978}{58}.}
  \ref{Hill}{Hill,M.J.M. {\it Trans.Camb.Phil.Soc.} {\bf 13} (1883) 36.}
  \ref{Cayley}{Cayley,A. {\it Quart.Math.J.} {\bf 7} (1866) 304.}
  \ref{Seade}{Seade,J.A. {\it Anal.Inst.Mat.Univ.Nac.Aut\'on
  M\'exico} {\bf 21} (1981) 129.}
  \ref{CM}{Cisneros--Molina,J.L. {\it Geom.Dedicata} {\bf84} (2001)
  \ref{Goette1}{Goette,S. \jram {526} {2000} 181.}
  207.}
  \ref{NandO}{Nash,C. and O'Connor,D--J, \jmp {36}{1995}{1462}.}
  \ref{Dows}{Dowker,J.S. \aop{71}{1972}{577}; Dowker,J.S. and Pettengill,D.F.
  \jpa{7}{1974}{1527}; J.S.Dowker in {\it Quantum Gravity}, edited by
  S. C. Christensen (Hilger,Bristol,1984)}
  \ref{Jadhav2}{Jadhav,S.P. \prD{43}{1991}{2656}.}
  \ref{Dow11}{Dowker,J.S. \cqg{21}{2004}4247.}
  \ref{Dow12}{Dowker,J.S. \cqg{21}{2004}4977.}
  \ref{Dow13}{Dowker,J.S. \jpa{38}{2005}1049.}
  \ref{Zagier}{Zagier,D. \ma{202}{1973}{149}}
  \ref{RandG}{Rademacher, H. and Grosswald,E. {\it Dedekind Sums},
  (Carus, MAA, 1972).}
  \ref{Berndt7}{Berndt,B, \aim{23}{1977}{285}.}
  \ref{HKMM}{Harvey,J.A., Kutasov,D., Martinec,E.J. and Moore,G.
  {\it Localised Tachyons and RG Flows}, hep-th/0111154.}
  \ref{Beck}{Beck,M., {\it Dedekind Cotangent Sums}, {\it Acta Arithmetica}
  {\bf 109} (2003) 109-139 ; math.NT/0112077.}
  \ref{McInnes}{McInnes,B. {\it APS instability and the topology of the brane
  world}, hep-th/0401035.}
  \ref{BHS}{Brevik,I, Herikstad,R. and Skriudalen,S. {\it Entropy Bound for the
  TM Electromagnetic Field in the Half Einstein Universe}; hep-th/0508123.}
  \ref{BandO}{Brevik,I. and Owe,C.  \prD{55}{4689}{1997}.}
  \ref{Kenn}{Kennedy,G. Thesis. University of Manchester 1978.}
  \ref{KandU}{Kennedy,G. and Unwin S. \jpa{12}{L253}{1980}.}
  \ref{BandO1}{Bayin,S.S.and Ozcan,M.
  \prD{48}{2806}{1993}; \prD{49}{5313}{1994}.}
  \ref{Chang}{Chang, P. Thesis. University of Manchester 1993.}
  \ref{Barnesa}{Barnes,E.W. {\it Trans. Camb. Phil. Soc.} {\bf 19} (1903) 374.}
  \ref{Barnesb}{Barnes,E.W. {\it Trans. Camb. Phil. Soc.}
  {\bf 19} (1903) 426.}
  \ref{Stanley1}{Stanley,R.P. \joa {49Hilf}{1977}{134}.}
  \ref{Stanley}{Stanley,R.P. \bams {1}{1979}{475}.}
  \ref{Hurley}{Hurley,A.C. \pcps {47}{1951}{51}.}
  \ref{IandK}{Iwasaki,I. and Katase,K. {\it Proc.Japan Acad. Ser} {\bf A55}
  (1979) 141.}
  \ref{IandT}{Ikeda,A. and Taniguchi,Y. {\it Osaka J. Math.} {\bf 15} (1978)
  515.}
  \ref{GandM}{Gallot,S. and Meyer,D. \jmpa{54}{1975}{259}.}
  \ref{Flatto}{Flatto,L. {\it Enseign. Math.} {\bf 24} (1978) 237.}
  \ref{OandT}{Orlik,P and Terao,H. {\it Arrangements of Hyperplanes},
  Grundlehren der Math. Wiss. {\bf 300}, (Springer--Verlag, 1992).}
  \ref{Shepler}{Shepler,A.V. \joa{220}{1999}{314}.}
  \ref{SandT}{Solomon,L. and Terao,H. \cmh {73}{1998}{237}.}
  \ref{Vass}{Vassilevich, D.V. \plb {348}{1995}39.}
  \ref{Vass2}{Vassilevich, D.V. \jmp {36}{1995}3174.}
  \ref{CandH}{Camporesi,R. and Higuchi,A. {\it J.Geom. and Physics}
  {\bf 15} (1994) 57.}
  \ref{Solomon2}{Solomon,L. \tams{113}{1964}{274}.}
  \ref{Solomon}{Solomon,L. {\it Nagoya Math. J.} {\bf 22} (1963) 57.}
  \ref{Obukhov}{Obukhov,Yu.N. \pl{109B}{1982}{195}.}
  \ref{BGH}{Bernasconi,F., Graf,G.M. and Hasler,D. {\it The heat kernel
  expansion for the electromagnetic field in a cavity}; math-ph/0302035.}
  \ref{Baltes}{Baltes,H.P. \prA {6}{1972}{2252}.}
  \ref{BaandH}{Baltes.H.P and Hilf,E.R. {\it Spectra of Finite Systems}
  (Bibliographisches Institut, Mannheim, 1976).}
  \ref{Ray}{Ray,D.B. \aim{4}{1970}{109}.}
  \ref{Hirzebruch}{Hirzebruch,F. {\it Topological methods in algebraic
  geometry} (Springer-- Verlag,\break  Berlin, 1978). }
  \ref{BBG}{Bla\v{z}i\'c,N., Bokan,N. and Gilkey, P.B. {\it Ind.J.Pure and
  Appl.Math.} {\bf 23} (1992) 103.}
  \ref{WandWi}{Weck,N. and Witsch,K.J. {\it Math.Meth.Appl.Sci.} {\bf 17}
  (1994) 1017.}
  \ref{Norlund}{N\"orlund,N.E. \am{43}{1922}{121}.}
  \ref{Duff}{Duff,G.F.D. \aom{56}{1952}{115}.}
  \ref{DandS}{Duff,G.F.D. and Spencer,D.C. \aom{45}{1951}{128}.}
  \ref{BGM}{Berger,M., Gauduchon,P. and Mazet,E. {\it Lect.Notes.Math.}
  {\bf 194} (1971) 1. }
  \ref{Patodi}{Patodi,V.K. \jdg{5}{1971}{233}.}
  \ref{GandS}{G\"unther,P. and Schimming,R. \jdg{12}{1977}{599}.}
  \ref{MandS}{McKean,H.P. and Singer,I.M. \jdg{1}{1967}{43}.}
  \ref{Conner}{Conner,P.E. {\it Mem.Am.Math.Soc.} {\bf 20} (1956).}
  \ref{Gilkey2}{Gilkey,P.B. \aim {15}{1975}{334}.}
  \ref{MandP}{Moss,I.G. and Poletti,S.J. \plb{333}{1994}{326}.}
  \ref{BKD}{Bordag,M., Kirsten,K. and Dowker,J.S. \cmp{182}{1996}{371}.}
  \ref{RandO}{Rubin,M.A. and Ordonez,C. \jmp{25}{1984}{2888}.}
  \ref{BaandD}{Balian,R. and Duplantier,B. \aop {112}{1978}{165}.}
  \ref{Kennedy2}{Kennedy,G. \aop{138}{1982}{353}.}
  \ref{DandKi2}{Dowker,J.S. and Kirsten, K. {\it Analysis and Appl.}
 {\bf 3} (2005) 45.}
  \ref{Dow40}{Dowker,J.S. {\it p-form spectra and Casimir energy}
  hep-th/0510248.}
  \ref{BandHe}{Br\"uning,J. and Heintze,E. {\it Duke Math.J.} {\bf 51} (1984)
   959.}
  \ref{Dowl}{Dowker,J.S. {\it Functional determinants on M\"obius corners};
    Proceedings, `Quantum field theory under
    the influence of external conditions', 111-121,Leipzig 1995.}
  \ref{Dowqg}{Dowker,J.S. in {\it Quantum Gravity}, edited by
  S. C. Christensen (Hilger, Bristol, 1984).}
  \ref{Dowit}{Dowker,J.S. \jpa{11}{1978}{347}.}
  \ref{Kane}{Kane,R. {\it Reflection Groups and Invariant Theory} (Springer,
  New York, 2001).}
  \ref{Sturmfels}{Sturmfels,B. {\it Algorithms in Invariant Theory}
  (Springer, Vienna, 1993).}
  \ref{Bourbaki}{Bourbaki,N. {\it Groupes et Alg\`ebres de Lie}  Chap.III, IV
  (Hermann, Paris, 1968).}
  \ref{SandTy}{Schwarz,A.S. and Tyupkin, Yu.S. \np{242}{1984}{436}.}
  \ref{Reuter}{Reuter,M. \prD{37}{1988}{1456}.}
  \ref{EGH}{Eguchi,T. Gilkey,P.B. and Hanson,A.J. \prp{66}{1980}{213}.}
  \ref{DandCh}{Dowker,J.S. and Chang,Peter, \prD{46}{1992}{3458}.}
  \ref{APS}{Atiyah M., Patodi and Singer,I.\mpcps{77}{1975}{43}.}
  \ref{Donnelly}{Donnelly.H. {\it Indiana U. Math.J.} {\bf 27} (1978) 889.}
  \ref{Katase}{Katase,K. {\it Proc.Jap.Acad.} {\bf 57} (1981) 233.}
  \ref{Gilkey3}{Gilkey,P.B.\invm{76}{1984}{309}.}
  \ref{Degeratu}{Degeratu.A. {\it Eta--Invariants and Molien Series for
  Unimodular Groups}, Thesis MIT, 2001.}
  \ref{Seeley}{Seeley,R. \ijmp {A\bf18}{2003}{2197}.}
\end{putreferences}

\bye